\magnification=\magstep1
\input amstex
\documentstyle{amsppt}

\hsize=5in
\vsize=7.5in

\def\bbC{\Bbb C}
\def\bbZ{\Bbb Z}

\def\bbR{\Bbb R}

\def\calo{\Cal O}

\def\cK{\Cal K}
\def\cL{\Cal L}
\def\cM{\Cal M}

\def\cN{\Cal N}

\def\D{\operatorname{D}}

\def\Hom{\operatorname{Hom}}
\def\Ext{\operatorname{Ext}}

\def\isomo{\overset{\sim}\to{=}}
\def\orig{\underline{\text{o}}}
\def\Stab{\operatorname{Stab}}
\def\spane{\operatorname{span}}
\def\Star{\operatorname{Star}}
\def\pf{\hfill $\square$}

\topmatter

\title
	{Toric varieties and minimal complexes}
\endtitle

\author
	{P.Bressler and V.Lunts}
\endauthor
\address
	{Department of Mathematics, The Pennsylvania State University,
	University Park, PA 16802}
\endaddress
\email
	{bressler\@math.psu.edu}
\endemail
\address
	{Department of Mathematics, Indiana University,
	Bloomington, IN 47405}
\endaddress
\email
	{vlunts\@indiana.edu}
\endemail

\thanks
	{The second author was supported in part by NSF grant 
	DMS-9504522.}
\endthanks
\endtopmatter

\document

\head{\bf 1. \ Introduction}\endhead

Let $T = \left(\bbC^\ast\right)^n$ and $N:=\operatorname{Lie}_{\bbR}(T)
\simeq\bbR^n$. Let $A=\operatorname{Sym}^\bullet N^*$ be the graded
polynomial ring of functions on $N$ where linear functions have degree 2.
Let $BT$ be the classifying space for $T$. There exists a natural
(and well known) isomorphism of graded algebras
$$A\simeq H(BT).$$

Let $X$ be a ($T$-)toric variety. The above isomorphism suggests that the
equivariant geometry on $X$ can be described in terms of coherent
objects on $N$, or more precisely, on the fan $\Phi$ corresponding to $X$.

In fact, in [BL],ch.15 a certain natural "minimal" complex $\cK_{\Phi}$
of $A$-modules was defined on $\Phi$. In case $\Phi$ is a complete fan,
$\cK_{\Phi}$ has remarkable properties: it is acyclic, except in degree
$-n$ and $H^{-n}(\cK_{\Phi})\simeq IH_{T,c}(X)$ - the equivariant
intersection cohomology with compact supports of $X$ (see [BL]). In a
sense, the coherent complex $\cK_{\Phi}$ on $\Phi$ represents the
equivariant intersection cohomology complex $IC_T(X)$ on $X$.

Let $\pi:X \to Y$ be a proper morphism of toric varieties. Denote
also by $\pi:\Phi \to \Psi$ the subdivision of corresponding fans
($\Phi$ subdivides $\Psi$). By the equivariant decomposition theorem
([BL]) the direct image $\pi_*IC_T(X)$ is isomorphic to a direct sum
of simple equivariant perverse sheaves on $Y$. In this work we translate
this theorem into a statement about the semisimplicity of the
direct image $\pi_*\cK_{\Phi}$ (which we define first) (Theorem 3.3.1).

The authors would like to thank the MIT Math Department common room for
providing coffee and a stimulating working environment.

\head{\bf 2. \ Complexes associated to fans}\endhead

\subhead{2.1. \ Fans and toric varieties}\endsubhead
Let $T = \left(\bbC^\ast\right)^n$ be a complex torus of dimension $n$,
$\Lambda = \Hom (\bbC^\ast,T)\isomo\bbZ^n$, $N=\Lambda\otimes _{{\Bbb Z}}\bbR
\isomo\bbR^n$. The space $N\otimes _{{\Bbb R}}\bbC$ is naturally isomorhic to
the Lie algebra of $T$.

A {\it convex polyhedral cone in $N$} (or, simply a {\it cone})
is a closed subset of $N$ which is the collection
of solutions of a finite system of linear homogeneous equations and
inequalities.
The term "convex" means "strictly convex", i.e. we assume that a cone
does not contain a nonzero linear subspace.
A cone is called rational if the linear functionals in the systems
take rational values on $\Lambda$. Let $\orig$ denote the cone which
consists of the origin of $N$.

A {\it fan} $\Phi$ in $N$ is a (finite) collection of convex polyhedral
cones in $N$, such that, if $\sigma\in\Phi$ and $\tau$ is a face of $\sigma$,
then $\tau\in\Phi$. Also any two cones in $\Phi$ intersect along a
common face.
In particular $\orig$ belongs to every fan.

For a cone $\sigma\subset N$ we denote by $\langle\sigma\rangle$
the fan generated by $\sigma$ (i.e. consisting of $\sigma$ and its faces),
and by $\partial\sigma$ the fan which consists of the proper faces of
$\sigma$.

A fan $\Phi$ is called {\it rational} is every cone
of $\Phi$ is rational. A fan $\Phi$ is called {\it complete} if the union of
all cones of $\Phi$ is equal to $N$.

Let $A=\bbR[x_1,\dots,x_n]$ denote the ring of polynomial functions
on $N$ considered as a graded $\bbR$-algebra by $\deg(x_i)=2$. Let
$\frak{m}$ denote the maximal ideal of elements of positive degree.

For a cone $\sigma\subset N$ let $A_\sigma$ denote the graded ring of
polynomial functions on $\sigma$. An inclusion of cones
$\tau\subset\sigma$ induces the restriction homomorphism of rings $\CD
A_\sigma \to A_\tau\endCD$.  In particular there are canonical
restriction maps $\CD A \to A_\sigma\endCD$ for all cones $\sigma$.

An irreducible algebraic veriety $X$ is called {\it $T$-toric} (or, simply,
{\it toric}), if $T$ acts on $X$ and $X$ contains a dense $T$-orbit
isomorphic to $T$ (in which case the number of $T$-orbits is finite).

There is a natural one-to-one correspondence between normal toric
varieties and rational fans in $N$. We refer the reader to [KKMS-D] for
detail of this construction. In what follows we will denote by
$X_\Phi$ the normal toric variety which corresponds to the rational fan
$\Phi$ and by $\Phi _X$ the fan, corresponding to the toric variety $X$.
We will need the following facts and notations.

The orbits of $T$ on $X_\Phi$ are in one-to-one correspondence with the
cones of $\Phi$. Let $\calo_\sigma$ denote the orbit associated with the cone
$\sigma$, and let $\sigma_{\calo}$ denote the cone associated with the orbit
$\calo$. The subspace $\spane(\sigma_{\calo})\otimes _{{\Bbb R}}\bbC\subset
N\otimes _{{\Bbb R}}\bbC$ is identified with the Lie algebra of the
stabilizer
$\Stab(\calo)\subset T$ of $\calo$ and  $\dim _{{\Bbb C}}\calo = n -\dim
_{{\Bbb R}}\sigma$.

The ring $A$ is canonically isomorphic to the (graded) cohomology ring
$H(BT;\bbR)$ (where $BT$ is the classifying space of $T$). Under the
above correspondence the ring $A_\sigma$ is canonically isomorphic to
$H(B\Stab(\calo);\bbR)$.

A toric variety $X$ is {\it affine} (resp. {\it complete}),
 if $\Phi _X=\langle\sigma\rangle$ for some cone $\sigma$ (resp. the
union of cones in $\Phi _X$ is the whole space $N$).

Let $\Phi $, $\Psi $ be two fans in $N$. We say that $\Phi $ {\it maps}
to $\Psi $, if for every $\sigma \in \Phi $ there exists $\tau \in \Psi$,
such that $\sigma \subset \tau$. In case the fans are rational, $\Phi $
maps to $\Psi $ iff there exists a morphism of toric varieties
$X_{\Phi}\to X_{\Psi}$ (this morphism is unique). Moreover, this
morphism is {\it proper} iff
$$\bigcup _{\Sb{\sigma \in \Phi}\endSb}\sigma =
\bigcup _{\Sb{\tau \in \Psi}\endSb}\tau.
$$
We call $\Phi $ a {\it subfan} of $\Psi $ iff every cone in $\Phi $
belongs to $\Psi $. This corresponds to an {\it open embedding} of
corresponding toric varieties.

\subhead{2.2. \ Complexes on fans}\endsubhead
An "$A$-module" will mean "a graded $A$-module". Morphisms between
$A$-modules preserve the grading. The tensor product $\otimes $ means
$\otimes _{\bbR}$ unless specified otherwise.
Given an $A$-module $M$ we denote by $\overline{M}$
the graded vector space $M/\frak{m}M$.

Fix a fan $\Phi$ in $N$, not necessarily rational.

\demo
{\bf Definition 2.2.1}
A {\it complex $\cM^\bullet$ on the fan $\Phi$} is a complex of 
$A$-modules 
$$
\CD
	0 @>>> \cM^{-n} @>{d^{-n}}>> \cM^{-n+1} @>{d^{-n+1}}>>
	\dots @>{d^{-1}}>> \cM^0 @>>> 0
\endCD
$$
such that
\roster
\item for each $i$ the module $\cM^i$ has a direct sum decomposition
$$
       	\cM^{-i} = \bigoplus_{\Sb{\sigma\in\Phi} \\ 
	{\dim\sigma=i}\endSb}\cM_\sigma
$$
\item each summand $\cM_\sigma$ is a finitely generated 
$A_\sigma$-module;
\item the differential $d:\cM^{\bullet} @>>> \cM^{\bullet+1}$ restricts to a
morphism of modules
$$
	d: \cM_\sigma @>>> 
	\underset{\Sb{\tau\in\partial\sigma} \\ {\dim\tau=\dim\sigma-1}\endSb}
	\to{\bigoplus}
	\cM_\tau
$$
compatible with the restriction maps of rings $A_\sigma @>>> A_\tau$.
\endroster
\enddemo

Suppose that $\Psi$ is a subfan of $\Phi$.  Then, for a complex $\cM^\bullet$
on $\Phi$, the submodule
$\bigoplus_{\sigma\in\Psi}\cM_\sigma\subset\bigoplus_i\cM^{-i}$ is preserved
by the differential.

\demo
{\bf Definition 2.2.2}
The submodule $\bigoplus_{\sigma\in\Psi}\cM_\sigma$ endowed with the
restriction of the differential in $\cM^\bullet$ is called {\it the
restriction of $\cM^\bullet$ to $\Psi$} and will be denoted by
$\left.\cM^\bullet\right\vert_\Psi$.
\enddemo

\demo
{\bf Definition 2.2.3}
A complex $\cM^\bullet$ on the fan $\Phi$ is called {\it locally free}
if for every $\sigma\in\Phi$ the $A_\sigma$-module $\cM_\sigma$ is free.

\enddemo
\demo
{\bf Definition 2.2.4}
A complex $\cM^\bullet$ on the fan $\Phi$ is called {\it locally exact}
if for every $\tau\in\Phi$  of positive dimension the complex
$\left.\cM^\bullet\right\vert_{\langle\tau\rangle}$
is acyclic in degree $-\dim\tau+1$, i.e. the differential
$$
\CD
\cM_{\tau}
@>{d}>>
\ker\left(\left(\cM^\bullet\big\vert_{\langle\tau\rangle}\right)^{-dim\tau+1}
\right.
@>{d}>> 
\left.
\left(\cM^\bullet\big\vert_{\langle\tau\rangle}\right)^{-dim\tau+2}
\right)
\endCD
\tag2.1
$$
is a surjection.
\enddemo

Note that the last two definitions are local, i.e. stable under
restriction to a subfan.

\demo
{\bf Definition 2.2.5.} For a complex $\cM^\bullet$ on the fan $\Phi$
and a cone $\sigma \in \Phi$ we call $\cM_\sigma$ the $\sigma $-
{\it component} of $\cM^\bullet$. The {\it support} of $\cM^\bullet$
is the union of the cones $\sigma \in \Phi$, such that $\cM_\sigma \neq 0$.

\enddemo

\subhead{2.3. \ Minimal complexes}\endsubhead
Recall that $\orig$ denotes the origin in $N$ viewed as a zero-dimensional
cone. Thus $A_{\orig}\isomo\bbR$. Let $\bbR(n)$ denote the free 
$A_{\orig}$-module with a unique generator of degree $-n = -\dim_{\bbC}T$.
Recall that $\frak{m}$ denotes the ideal of elements of positive degree in
$A$.

\demo
{\bf Definition 2.3.1}
A complex $\cM^\bullet$ on the fan $\Phi$ is called {\it minimal}
if it satisfies the following:
\roster
\item $\cM_{\orig} = \bbR(n)$;
\item it is locally free and locally exact;
\item for every $\tau\in\Phi$ of dimension $\dim\tau>0$ the reduction of
the map (2.1) modulo $\frak{m}$ is an isomorphism.
\endroster
\enddemo

For any fan $\Phi$ one can construct a minimal complex on $\Phi$. Moreover,
one can show that any two minimal complexes are non-canonically isomorphic.
In what follows we will denote a minimal complex on $\Phi$ by $\cK_\Phi$
or simply by $\cK$. The definition of a minimal complex is local, i.e. the
restriction of a minimal complex to a subfan is also minimal.

The remarkable properties of minimal complexes are summarized in the following
theorem (Theorem 15.7 of \cite{BL}).
\proclaim
{Theorem 2.3.2}
Suppose that $\Phi$ is a rational fan in $N$ which is {\it either} complete
{\it or} generated by a single cone of dimension $n$. Then
\roster
\item the minimal complex $\cK_\Phi$ is acyclic except in the lowest degree
(which is equal to $-n$);
\item there is an isomorphism of graded $A$-modules
$H^{-n}\cK_\Phi\isomo IH_{T,c}(X_\Phi;\bbR)$;
\item for any cone $\sigma\in\Phi$ there is an isomorphism of
graded $A_\sigma$-modules $\cK_\sigma\isomo
H\left(IC(X_\Phi;\bbR)_{\calo_\sigma}\right)\otimes A_\sigma$.
\endroster
\endproclaim

The notations in 2.3.2 are as follows:
$X_\Phi$ is the normal toric variety which corresponds to the rational
fan $\Phi$; $IH_{T,c}(X_\Phi;\bbR)$ is
the $T$-equivariant intersection cohomology of $X_\Phi$ with compact supports
(\cite{BL},13.4);
$IC(X_\Phi;\bbR)_{\calo_\sigma}$ is the stalk of the intersection complex
of $X_\Phi$ at a point of the orbit $\calo_\sigma$ which corresponds to the 
cone $\sigma$ -- this is a complex of vector spaces with naturally graded
(by cohomological degree) finite dimensional cohomology.

\remark
{Remark 2.3.3} We know of no "elementary" proof of the above theorem.
In particular, we cannot prove the first statement for fans, which
are not rational.
 The proof of Theorem 2.3.2 depends on the 
following fact:
{\it
if $\Phi$ is a rational fan in $N$ which is either complete
or generated by a single cone of dimension $n$, then
$IH_{T,c}(X_\Phi;\bbR)$ is a free module over $A$.} 
The same proof shows that the conclusions of 2.3.2 hold
for {\it any} rational fan $\Phi$ such that $IH_{T,c}(X_\Phi;\bbR)$ is
a free $A$-module.
\endremark

\proclaim
{Proposition 2.3.4.} 
Suppose that $\Phi$ is a rational fan in $N$ generated by several
$n$-dimensional cones whose union is convex (therefore a convex cone).
Then $IH_{T,c}(X_\Phi;\bbR)$ is a free module over $A$.
\endproclaim

\demo
{Proof}
Let $\sigma$ denote the union of the $n$-dimensional cones of $\Phi$.
Then $\sigma$ is a convex cone in $N$ and $\Phi$ is a subdivision of
$\langle\sigma\rangle$, hence there is a unique morphism of toric varieties
$$
\pi : X_\Phi @>>> X_{\langle\sigma\rangle}
$$
which is $T$-equivariant and proper. Therefore the direct image functors
$$
\pi_*,\ \pi_! : \D_{T,c}(X_\Phi) @>>> \D_{T,c}(X_{\langle\sigma\rangle})
$$
are defined and $\pi_* = \pi_!$. Here $D_{T,c}$ denotes the bounded
constructible equivariant derived category of sheaves as defined in
\cite{BL}. 
 By the equivariant decomposition theorem
(\cite{BL}, 5.3) the object $\pi_!IC_T(X_\Phi)$ is isomorphic to
a direct sum of shifted simple equivariant perverse sheaves on
$X_{\langle\sigma\rangle}$. By \cite{BL}, 5.2, simple equivariant perverse
sheaves on $X_{\langle\sigma\rangle}$ are of the form
$IC_T(Z)$, where $Z$ is a closed $T$-invariant irreducible subvariety of
$X_{\langle\sigma\rangle}$. (Here we use the normality of
$X_{\langle\sigma\rangle}$ to conclude that the stabilizers of all points
in $X_{\langle\sigma\rangle}$ are connected, hence the only $T$ equivariant
local systems on $T$-invariant subvarieties are the trivial ones.)

Consider a closed irreducible $T$-invariant subvariety $Z$ of
$X_{\langle\sigma\rangle}$. Let $T_0$ denote the stabilizer of a general
point of $Z$, and let $T_1=T/T_0$. Then, $Z$ is a $T_1$-toric variety.
Since $X_{\langle\sigma\rangle}$ is affine with a unique $T$-fixed point
$Z$ is also affine with a unique $T_1$-fixed point.

Therefore, by \cite{BL} 14.3 (ii'), $IH_{T_1,c}(Z)$ is a free module over
$H(BT_1;\bbR)$. Since $A=H(BT;\bbR)\isomo H(BT_0;\bbR)\otimes
H(BT_1;\bbR)$ and $T_0$ acts trivially on $Z$, it follows that
$IH_{T_,c}(Z)$ is a free module over $A$.

Hence, $IH_{T,c}(X_\Phi;\bbR)\isomo H_{T,c}(X_{\langle\sigma\rangle};
\pi_!IC_T(X_\Phi))$ is a free module over $A$.  \pf
\enddemo

\proclaim
{Corollary 2.3.5}
Suppose that $\Phi$ is a rational fan in $N$ generated by several
$n$-dimensional cones whose union is convex (therefore a convex cone).
Then the conclusions of Theorem 2.3.2 hold for $\Phi$.
\endproclaim

\demo
{Proof}
Follows from Proposition 2.3.4 by Remark 2.3.3.  \pf
\enddemo

\subhead{2.4. \ Shifted minimal complexes}\endsubhead
Let $\Phi$ be a fan in $N$ (not necessarily rational) and let $\sigma$
be a cone in $\Phi$. We have the following generalization of the
minimal complex $\cK_\Phi$.
\demo
{\bf Definition 2.4.1}
A {\it shifted (by $k$) minimal complex based on $\sigma$}, denoted by
\linebreak
$\cK_\Phi[\sigma](k)$, is a complex on $\Phi$ which satisfies the following
conditions:
\roster
\item $\cK_\Phi[\sigma](k)$ is locally free and locally exact;

\item $\cK_\Phi[\sigma](k)^{-i} = 0$ for $i<\dim\sigma$;

\item $\cK_\Phi[\sigma](k)^{-\dim\sigma} = (\cK_\Phi[\sigma](k))_\sigma =
A_\sigma(n-\dim\sigma+k)$ (where $A_\sigma(n-\dim\sigma+k)$ is the free graded
$A_\sigma$-module on one generator of degree $n-\dim\sigma+k$);

\item for every $\tau\in\Phi$ of dimension $\dim\tau>\dim\sigma$ the
reduction of the map 2.1 (with
$\cM^\bullet=\cK_\Phi[\sigma](k)$) modulo $\frak{m}$ is an isomorphism.
\endroster
\enddemo

\remark
{Remark 2.4.2}
\roster
\item It is easy to show that $\sigma$ and $k$ determine $\cK_\Phi[\sigma](k)$
uniquely up to a non-canonical isomorphism.
\item $\cK_\Phi = \cK_\Phi[\orig](0)$.
\endroster
\endremark

Recall that, for a cone $\sigma$ of $\Phi$ the {\it star of $\sigma$} is
defined as
$$
	\Star(\sigma)\overset{\text{def}}\to{=}\left\lbrace\tau\in
	\Phi\left.\right\vert
	\sigma\in\langle\tau\rangle\right\rbrace\ .
$$
As is easy to see, the shifted minimal complex $\cK_\Phi[\sigma](k)$ is
supported on $\Star(\sigma)$.

\proclaim
{Lemma 2.4.3}
Suppose that $\Phi$ is a rational fan which is {\it either} complete {\it or}
is generated by a single cone of dimension $n$, $\sigma$ is a cone in $\Phi$,
and $k$ is an integer. Then the shifted minimal complex $\cK_\Phi[\sigma](k)$
is acyclic except in the lowest degree (which is equal to $-n$).
\endproclaim

\demo
{Proof}
Let $p:N @>>> N^{\prime }:=N/\spane(\sigma)$ denote the quotient map.
The image $\Psi$ of
$\Star(\sigma)$ under $p$ is a rational fan in $N^{\prime}$ which is
complete (respectively generated by a single cone of dimension
$n-\dim\sigma$) if $\Phi$ is complete (respectively generated by a single cone
of dimension $n$).

Let $A'$ denote the ring of polynomial functions on $N^{\prime}$.
The map $p$ induces on $A$ a structure of an algebra over $A'$. Let
$\orig ^{\prime}$ be the origin in $N^{\prime }$. 
 It is
easy to check that
$$
\cK_\Phi[\sigma](k) = \left(A\otimes_{A'}\cK_\Psi[\orig ^{\prime}](k)
\right)[\dim\sigma] \ 
$$
as complexes of $A$-modules.
By Theorem 2.3.2 the complex $\cK_\Psi[\orig ^{\prime}](k)$ is acyclic
except in degree $-n+\dim\sigma$. Hence $\cK_\Phi[\sigma](k)$ is acyclic
except in degree $-n$.  \pf
\enddemo

\subhead{2.5. \ A Criterion for semi-simplicity}\endsubhead

\proclaim
{Lemma 2.5.1}
Suppose that $\Phi$ is a fan in $N$ (not necessarily rational) and
$\cM^\bullet$ is a locally free, locally exact complex on $\Phi$,
 such that
$\cM^p=0$ for $p>-i$. 
Let $\tau \in \Phi$ be of dimension $i$. Then $\cM^\bullet$ contains
a direct summand isomorphic to $\cK_{\Phi}^{\bullet}[\tau ](i-n)
\otimes\overline{\cM}_{\tau}$.

\endproclaim

\demo
{Proof} For a complex $\cL ^{\bullet}$ denote by $\cL ^{>k}$ the
truncated complex
$$0 @>>> \cL^{k+1} @>>> \cL^{k+2} @>>> \cdots \ . $$

We will construct a splitting
$$\cM^{\bullet}\simeq (\cK_{\Phi}^{\bullet}[\tau ](i-n)
\otimes\overline{\cM}_{\tau})\oplus \cN^{\bullet}$$
by induction on the dimension of cones working from right to left.
Since $\cK_{\Phi}^{-i}[\tau ](i-n)
\otimes\overline{\cM}_{\tau}\simeq \overline{\cM}_{\tau}$, the
decomposition
$$\cM^{-i}=\cM_{\tau}\oplus \bigoplus _{\Sb{dim(\xi)=i} \\
{\xi \neq \tau}\endSb}\cM_{\xi}
$$
determines the desired decomposition
$$\cM^{\geq -i}=(\cK_{\Phi}^{\bullet}[\tau ](i-n)
\otimes\overline{\cM}_{\tau})^{\geq -i}\oplus \cN^{\geq -i}.
$$

Assume that we have constucted a splitting
$$\cM^{> -j}=(\cK_{\Phi}^{\bullet}[\tau ](i-n)
\otimes\overline{\cM}_{\tau})^{> -j}\oplus \cN^{> -j}
$$
for some $j>i$. Let $\sigma \in \Phi$ be of dimension $j$. It suffices
to show that we can extend the splitting 
to the subfan $\langle\sigma\rangle$. The following claim ensures that
we can do so.

\proclaim
{Claim 2.5.2}
Let
$$d:L \twoheadrightarrow Z
$$
be a surjection of finitely generated (graded) $A_{\sigma }$-modules,
where the
module $L$ is free. Given a decomposition $Z=Z_K\oplus Z_N$,
there exists a decomposition $L=L_K\oplus L_N$ with the properties

a) $d(L_K)=Z_K,\ \ d(L_N)=Z_N$

b) $d$ induces an isomorphism of residues $\overline{d}:
\overline{L}_K @>{\sim}>>\overline{Z}_K$.

\endproclaim

\demo
{Proof of claim} Choose a subspace $\tilde{Z}_K\subset Z_K$, which
maps isomorphically onto $\overline{Z}_K$ under he residue map.
Choose a subspace  $S_K\subset L$ so that $d:S_K @>{\sim}>>
\tilde{Z}_K$. Since $S_K\cap {\frak m} L=0$, we can choose a subspace
$S_N\subset L$, such that $S_K\cap S_N=0$ and $S_K\oplus S_N$ generates
$L$ freely. Subtracting, if necessary, elements of $S_K$ from elements
of $S_N$ we may assume that $d(S_N)\subset N$. Thus we may take
$L_K:=A_{\sigma }S_K$, $L_N:=A_{\sigma }S_N$. This proves the claim and
the lemma. \pf

\enddemo
\enddemo

\proclaim
{Proposition 2.5.3}
Suppose that $\Phi$ is a fan in $N$ (not necessarily rational) and
$\cM^\bullet$ is a locally free, locally exact complex on $\Phi$.
Then $\cM^\bullet$ is noncanonically
isomorphic to a direct sum of shifted minimal complexes.
\endproclaim

\demo
{Proof}
The proposition follows by induction on $i=\min\left\lbrace p\vert
\cM^p\neq 0\right\rbrace$ and the number of cones of
dimension $i$ using Lemma 2.5.1.  \pf
\enddemo

\proclaim
{Corollary 2.5.4}
Suppose that $\Phi$ is a rational fan which is {\it either} complete {\it or}
is generated by a single cone of dimension $n$ and
$\cM^\bullet$ is a locally free, locally exact complex on $\Phi$.
Then $\cM^\bullet$ is acyclic except in the lowest degree (which is equal to
$-n$).
\endproclaim

\demo
{Proof}
By Proposition 2.5.3 $\cM^\bullet$ is (isomorphic to) a direct sum of
shifted minimal complexes. The conclusion now follows from Lemma
2.4.3.  \pf
\enddemo

\head{\bf 3. \ Direct image complexes and a decomposition theorem}\endhead

\noindent{\bf 3.0.} 
In this section we consider normal toric  varieties $X$ and $Y$
with corresponding rational fans $\Phi$ and $\Psi$ in $N$ and a proper
morphism of toric varieties
$$
\pi : X @>>> Y\ .
$$
Thus, $\Phi$ is a subdivision of $\Psi$. Put $\cK_X^\bullet=\cK_\Phi^\bullet$
and  $\cK_Y^\bullet=\cK_\Psi^\bullet$

By the equivariant decomposition theorem (\cite{BL}, 5.3) the object
$\pi_*IC_T(X)$ is (noncanonically isomorphic to) a direct sum of shifted
objects $IC_T(Z)$ for various closed $T$-invariant irreducible
subvarieties of $Y$
(see the proof of Proposition 2.3.4).

We will translate this into the language of minimal complexes.
Namely, we will
\roster
\item construct a direct image complex $\pi _*\cK^\bullet$ on
$\Psi$: we define $\pi _*\cK^\bullet$ as a certain subcomplex of
$\cK_X^\bullet$; 
\item show that the inclusion
$\pi _*\cK^\bullet\hookrightarrow\cK_X^\bullet$ is
a quasiisomorphism;
\item show that $\pi _*\cK^\bullet$ is (noncanonically isomorphic to) a
direct sum of shifted minimal complexes on $\Psi$;
\item in particular $\cK_Y^\bullet$ is a direct summand in
$\pi _*\cK^\bullet$.

\endroster

\subhead{3.1. \ Construction of the direct image complexes}\endsubhead
Suppose for a moment that $\Phi$ and $\Psi$ are fans in $N$ (not
necessarily rational) and $\Phi$ maps to $\Psi$. Call this map of fans
$\pi: \Phi @>>> \Psi$. Let $\cM^\bullet$ be a complex on $\Phi$. Let
us define its direct image $\pi _*\cM^\bullet$ on $\Psi$. It will be
defined as a subcomplex of $\cM^\bullet$. The terms of
$\pi _*\cM^\bullet$ are defined inductively as follows.

Suppose that
$(\pi _*\cM)^{>-i}$ has been defined (as a subcomplex of $\cM^{>-i}$).

For $\sigma\in\Psi$ let $\pi^{-1}(\sigma)$ denote the collection of cones
$\tau$ of $\Phi$ of the same dimension as $\sigma$ which
are contained in  $\sigma$.

Consider a cone $\sigma$ of dimension $i$. The requirement that
the square
$$
\CD
(\pi_*\cM)_\sigma @>{\left.d\right\vert_{(\pi_*\cM)_\sigma}}>>
(\pi_*\cM)^{-i+1} \\
@VVV @VVV \\
\displaystyle\bigoplus_{\tau\in\pi^{-1}(\sigma)}\cM_\tau @>{d}>>
\cM^{-i+1}
\endCD
$$
is Cartesian determines $\pi_*\cM_\sigma$  and $\left.d\right\vert_{
(\pi_*\cM)_\sigma}$
uniquely. Let
$$
(\pi_*\cM)^{-i}= \bigoplus_{\Sb{\sigma\in\Psi} \\ {\dim\sigma=i}\endSb}
(\pi_*\cM)_\sigma,
\ \ \ d=\sum_\sigma\left.d\right\vert_{(\pi_*\cM)_\sigma}\ .
$$
This defines $\pi_*\cM^\bullet$ as a complex on the fan $\Psi$.

Note, that the construction of $\pi_*\cM$ is local, i.e. compatible
with the restriction to a subcomplex $\Psi ^{\prime}\subset\Psi$
(and the corresponding subcomplex $\Phi^{\prime}=\Phi\cap\Psi^{\prime}$).

In the remainder of this section we will investigate the properties of
the direct image of the minimal complex under a proper morphism
of normal toric varieties.

\subhead{3.2. \ Properties of $\pi_*\cK_X$}\endsubhead
\medskip
We are back to notations of 3.0 above.

\proclaim
{Lemma 3.2.1}
Suppose that $\pi:X @>>> Y$ is a proper morphism of normal toric varieties.
Then the complex $\pi_*\cK_X^\bullet$ is locally exact.
\endproclaim

\demo
{Proof}
Fix a cone $\sigma\in\Psi$. To check local exactness of
$\pi_*\cK_X^\bullet$ at $\sigma$ we may assume that
$\Psi=\langle\sigma\rangle$. Clearly
$\Phi$ satisfies the hypotheses of Corollary
2.3.5 (with $N$ replaced by $\spane(\sigma)$ and $n$ replaced by
$\dim\sigma$) and therefore the minimal complex
 $\cK_X^\bullet$ is acyclic except in degree
$-\dim\sigma$.

Let $i=\dim\sigma$.
By construction of $\pi_*\cK_X^\bullet$, the square
$$
\CD
(\pi_*\cK^\bullet)_{\sigma} 
@>{d}>> 
(\pi_*\cK^\bullet)^{-i+1} \\
@VVV @VVV \\
\cK_X^{-i}
 @>{d}>>
\cK_X^{-i+1}
\endCD
$$
is Cartesian. Since $\cK^\bullet_X$
is acyclic in degree $-i+1$, this implies that
$\pi_*\cK^\bullet_X$ is acyclic in
degree $-i+1$.

As this happens for every cone $\sigma\in\Psi$, the complex
$\pi_*\cK_X^\bullet$
is locally exact.  \pf
\enddemo

\proclaim
{Proposition 3.2.2}
Suppose that $\pi:X @>>> Y$ is a proper morphism of normal toric varieties.
Then
\roster
\item the complex $\pi_*\cK_X^\bullet$ is locally free;
\item the inclusion $\pi_*\cK_X^\bullet\hookrightarrow\cK_X^\bullet$ is
a quasiisomorphism.
\endroster
\endproclaim

\demo
{Proof}
By induction on dimension of the torus we may assume that the conclusions
hold for any proper morphism of normal toric varieties of dimension smaller
than $n=\dim_{\bbR}N$.

\proclaim
{Case 1} The fan $\Psi$ is generated by a single cone $\sigma$
of dimension $i$, i.e.
$\Psi=\langle\sigma\rangle$.
\endproclaim

By induction on the dimension of the cone (the case of degree zero being
obvious since $A_{\orig}=\bbR$) we may assume that $(\pi_*\cK_X^\bullet)
_\tau$ is free over
$A_\tau$ for all $\tau \in \Psi$ of dimension smaller than $i$.

Denote by $\partial\Psi$ the subfan of $\Psi$ consisting of proper
faces of $\sigma$, and by $\partial\Phi$ its preimage in $\Phi$.
Consider the commutative diagram with rows short exact sequences
of complexes (which serves as the
definition of the complex $\cM^\bullet$):
$$
\CD
0 @>>> \pi_*\cK_X^\bullet\vert_{\partial\Psi} @>>>
\pi_*\cK_X^\bullet @>>>
(\pi_*\cK^\bullet)_\sigma[i] @>>> 0 \\
& & @VVV @VVV @VVV \\
0 @>>> \cK_X^\bullet\vert_{\partial\Phi} @>>>
\cK_X^\bullet @>>>
\cM^\bullet @>>> 0
\endCD
\tag3.1
$$

\proclaim
{Claim 3.2.3}
\roster
\item
The inclusion $\pi_*\cK_X^\bullet\vert_{\partial\Psi}\hookrightarrow
\cK_X^\bullet\vert_{\partial\Phi}$ (the left vertical
arrow in (3.1)) is a quasiisomorphism.
\item The complexes $\pi_*\cK_X^\bullet \vert_{\partial\Psi}$ and
$\cK_X^\bullet \vert_{\partial\Phi}$ are acyclic
except in the lowest degree $-i+1$.
\endroster
\endproclaim

\demo
{Proof}
Replacing $N$ by $\spane(\sigma)$ is necessary we may assume that
$\dim\sigma=n$. Choose an integral vector $v$ in the interior of $\sigma$,
and consider the projection $N @>>> N/\spane(v)\overset{\text{def}}\to{=}N'$.
Then, $\partial\Psi$ and $\partial\Phi$ project (isomorphically)
onto complete
rational fans in $N'$ which we denote $\Psi'$ and $\Phi'$ respectively.
Since $\Phi'$ is a subdivision of $\Psi'$ there is a proper morphism
of corresponding toric varieties
$$
	\pi' : X' @>>> Y'\ .
$$
Let $A'$ denote the ring of polynomial functions on $N'$, naturally
included in $A$. Then, clearly there are isomorphisms of complexes of
$A'$-modules $\cK_X^\bullet \vert_{\partial\Phi}\isomo
\cK_{X'}^\bullet$ and $\pi_*\cK_X^\bullet \vert_{\partial\Psi}
\isomo\pi'_*\cK_{X'}^\bullet$
compatible with inclusions of the latter into the former.

The first assertion now follows by induction on the dimension of the torus.
The second assertion follows from Theorem 2.3.2 applied to $\Phi'$. \pf
\enddemo

\proclaim
{Claim 3.2.4}
The inclusion $\pi_*\cK_X^\bullet\hookrightarrow\cK_X^\bullet$
(the middle vertical arrow in (3.1)) is a quasiisomorphism.
\endproclaim

\demo
{Proof}
By the previous claim the complex $\pi_*\cK_X^\bullet \vert
_{\partial\Psi}$ is acyclic, except in the lowest degree $-i+1$.
By Lemma 3.2.1 the complex $\pi_*\cK_X^\bullet$ is locally exact.
Hence $\pi_*\cK_X^\bullet$ is acyclic, except in the lowest
degree $-i$. But it follows immediately from the construction of
$\pi_*\cK_X^\bullet$ that
$H^{-i}\pi_*\cK^\bullet\isomo
H^{-i}\cK_X^\bullet$.
\pf
\enddemo

It follows from Claims 3.2.3 and 3.2.4 that
the map $(\pi_*\cK_K^\bullet)_\sigma[i] @>>> \cM^\bullet$
(the right vertical map in
(3.1)) is a quasiisomorphism, in particular the complex
$\cM^\bullet$ is acyclic except in the lowest degree equal to $-i$.

Note that $\cM^{-j} = \bigoplus\cM_\tau$, where $\dim\tau=j$ and
$\cM_\tau$ is a free $A_\tau$-module. Therefore
$\Ext^p_{A_\sigma}(\cM^{-j},A_\sigma) = 0$ for $p\neq i-j$.
The standard argument shows that
$H^{-i}\cM^\bullet$ satisfies
$\Ext^p_{A_\sigma}(H^{-i}\cM^\bullet,A_\sigma) = 0$ for $p\neq 0$.
This implies that $(\pi_*\cK_X^\bullet)_\sigma=H^{-i}\cM^{\bullet}$
is free over $A_\sigma$. \pf

\proclaim
{Case 2} $\Psi$ is general.
\endproclaim

In this case we still know the first assertion of the proposition,
since it is a local property and the proof in Case 1 applies.
It remains to prove the second assertion.

For a cone $\tau\in\Psi$ let $\Phi\cap\tau:=\{\xi\in\Phi \vert
\xi\subset\tau\}$ be the corresponding subfan of $\Phi$. Denote again
by $\pi:\Phi\cap\tau\to\langle\tau\rangle$ the corresponding subdivision
of fans. The inclusion $\pi_*\cK_X^\bullet\hookrightarrow\cK_X^\bullet$
restricts to the inclusion $\pi_*\cK_{\Phi\cap\tau}^\bullet
\hookrightarrow\cK_{\Phi\cap\tau}^\bullet$.

Consider the Cech resolution of the pair of
complexes $\pi_*\cK_X^\bullet\hookrightarrow\cK_X^\bullet$
$$ 
\CD
\cdots @>>> \bigoplus_{\Sb{\tau,\xi\in\Psi}\endSb}
\pi_*\cK_{\Phi\cap\tau\cap\xi}^\bullet @>>>
\bigoplus_{\Sb{\tau\in\Psi}\endSb}\pi_*\cK_{\Phi\cap\tau}^\bullet
@>>> \pi_*\cK_X^\bullet \\
@VVV @VVV @VVV @VVV\\
\cdots @>>> \bigoplus_{\Sb{\tau,\xi\in\Psi}\endSb}
\cK_{\Phi\cap\tau\cap\xi}^\bullet @>>>
\bigoplus_{\Sb{\tau\in\Psi}\endSb}\cK_{\Phi\cap\tau}^\bullet
@>>> \cK_X^\bullet 
\endCD
$$   

The two rows are exact and each vertical arrow, except the rightmost,
is a quasiisomorphism by the proof in Case 1 above. Hence the
rightmost arrow is also a quasiisomorphism. This completes the proof of
proposition 3.2.2. \pf

\enddemo
\proclaim
{Corollary 3.2.5.} Suppose that $\pi_*X @>>> Y$ is a proper morphism
of normal toric varieties. Then the direct image complex $\pi_*\cK_X^
\bullet$
is a direct sum of shifted minimal complexes. $\cK_Y^\bullet$ is a
direct summand of $\pi_*\cK_X^\bullet$ (in fact the only one
based on $\orig$).
\endproclaim

\demo
{Proof} The first assertion is a direct consequence of Lemma 3.2.1,
Proposition 3.2.2.(1), and of Proposition 2.5.3. The second assertion
is clear by the construction of $\pi_*\cK_X^\bullet$. \pf
\enddemo

\subhead{3.3. \ A decomposition theorem for minimal complexes} \endsubhead
\medskip

Let us summarize our main results. 

\proclaim
{Theorem 3.3.1}
Suppose that $\pi : X @>>> Y$ is a proper morphism of normal toric varieties.
Then the direct image complex $\pi_*\cK^\bullet$ (see (3.1)) has
the following properties:
\roster
\item $\pi_*\cK_X^\bullet$ is locally exact and locally free;
\item the inclusion $\pi_*\cK_X^\bullet\hookrightarrow\cK^\bullet_X$ is
a quasiisomorphism;
\item $\pi_*\cK_X^\bullet$ is a direct sum of shifted minimal complexes on
$\Psi_Y$;
\item $\cK^\bullet_Y$ is a direct summand of $\pi_*\cK_X^\bullet$
(in fact the only one based on $\orig$).
\endroster
\endproclaim

\demo
{Proof}
The first two assertions follow from Lemma 3.2.1 and
Proposition 3.2.2. The last two assertions are taken from Corollary
3.2.5.  \pf
\enddemo

\enddocument